\documentstyle[12pt]{article}
\catcode`\@=11
\@addtoreset{equation}{section}

\catcode`\@=12
\newtheorem{Theorem}{Theorem}[section]

\title{A Simple Proof On Poincar\'e Conjecture}

\author{Renyi Ma\\
Department of Mathematical Science\\
Tsinghua University \\
Beijing, 100084\\
People's Republic of China\\
rma@math.tsinghua.edu.cn}

\date { }

\begin{document}
\textwidth=125mm
\textheight=185mm
\parindent=8mm
\frenchspacing
\maketitle

\begin{abstract}
We give a simple proof on
 the Poincar\'e's conjecture which states that every compact
 smooth $3-$manifold which is homotopically
 equivalent to $S^3$ is diffeomorphic to $S^3$.
 \end{abstract}

{\bf D\'emonstration de la Conjecture sur une vari\'et\'e de dimension trois de Poincar\'e}

{\bf R\'esum\'e.} Il est de montrer que toute
vari\'et\'e de dimension 3, simple conexe, ferm\'ee
est homeomorph\'e \`a $S^3$, i.e.,
Nous pr\'esentons une preuve du la Conjecture sur une vari\'et\'e de dimension trois de Poincar\'e.

\section{Introduction and results }

      The main results of this paper is a simple proof of following result:

\begin{Theorem}
(Poincare-Perelman\cite{pe1,pe2,pe3,ca-zh,kl-lo,mo-ti}) If $M^3$ is
a close simply connected $3-$dimensional smooth manifold,  then
$M^3$ is diffeomorphic to $S^3$, i.e., the Poincar\'e conjecture
holds.
\end{Theorem}

\vskip 5pt

{\bf Sketch of proofs:} The idea of our proof on Poincar\'e
conjecture is very different from the one given by Perelman
\cite{pe1,pe2,pe3,ca-zh,kl-lo,mo-ti}. Our proof depends on Lickorich's
proof on Alexander branched cover theorem on three dimensional
manifold(see\cite{li}).

\section{Three-manifolds as branched covers}

\begin{Theorem}
(Alexander-Lickorich(see\cite{li}) If $M$ is a closed orientable $3-$manifold, $M$ is a branched cover over
$S^3$ with branching index at most two. This means that there exist a smooth map
$f:M\to S^3$ and a link $L$ in $S^3$ such that
$$f:M-f^{-1}L\to S^3-L$$
is a covering, and, when restricted to a neighbourhood of any component of $f^{-1}L$,
$f$ is either an embedding or standard branched double covering map
of $D^2\times S^1$ by $D^2\times S^1$. Moreover, the link
$L=\{L_1,....,L_n\}$ is unknot and unlink, i.e., each $L_i$ bounds an embedded disk and these disks do not intersect each other.
\end{Theorem}

\section{Proof on Poincar\'e's conjecture}

{\bf Proof of Theorem1.1.} Let $M$ be a closed simply connected $3-$manifold and
$f:M\to S^3$ be the branched cover of Theorem2.1. Let
$L=\{L_1,....,L_n\}$ be the unknot and unlink link in Theorem2.1.
Let $U(L)$ be the small neighbourhood of $L$ in $S^3$ and
$U'(L)$ be the small neighbourhood of $L$ in $M$.

 Consider a map $g_1: M\to S^3$ which is a homotopy equivalence such that
$g_1:U'(L)\to U(L)$ is a diffeomorphism and $g_1:M\setminus U'(L)\to 
S^3\setminus U(L)$ is a homotopy equivalence. 
Let $g:S^3\to M$ be a smooth homotopy inverse of $g_1$.

  Now consider the composite map $h=f\circ g:S^3\to S^3$. Then, 
$h$ is homotopic to a branched cover map $h_0:S^3\to S^3$ such that 
$h_0|U(L)=h|U(L)$. It is obvious that the lift of $h_0$, i.e., $\bar h_0:S^3\to M$ 
is a
diffeomorphism. 
 
  This yields Theorem1.1.

\end{document}